\documentclass[12pt]{amsart}
\usepackage{epsfig}
\usepackage{amsmath}
\usepackage{amsfonts}
\usepackage{amssymb}
\usepackage{epigraph, fancyhdr}

\def\J{\mathcal J}

\def\O{\mathcal O}

\def\1{\mathbf 1}
\def\M{{\overline{\mathcal M}}}

\def\ZZ{\mathbb Z}
\def\CC{\mathbb C}
\def\RR{\mathbb R}

\def\Res{\operatorname{Res}}

\def\hat{\widehat}
\def\tilde{\widetilde}
\def\p{\partial}
\def\a{\alpha}
\def\b{\beta}

\def\la{\lambda}
\def\gL{\Lambda}

\def\m{{\mathbf m}}

\def\Eu{\operatorname{Euler}}

\def\tr{\operatorname{tr}}
\def\Lie{\operatorname{Lie}}

\renewcommand{\Delta}{\triangle}

\title[Toric $q$-hypergeometric functions]
      {Permutation-equivariant \\ quantum K-theory V. \\
     Toric $q$-hypergeometric functions}

\author[A. Givental]{Alexander GIVENTAL}
\thanks{This material is based upon work supported by the National 
Science Foundation under Grant DMS-1007164, and by the IBS Center for Geometry 
and Physics, POSTECH, Korea.} 

\date{July 20, 2015}

\begin{document}

\begin{abstract}

We first retell in the K-theoretic context the heuristics of $S^1$-equivariant Floer theory on loop spaces which gives rise to $D_q$-module structures, and in the case of toric manifolds, vector bundles, or super-bundles to their explicit $q$-hypergeometric solutions. Then, using the fixed point localization technique developed in Parts II--IV, we prove that these $q$-hypergeometric solutions represent K-theoretic Gromov-Witten invariants.         

\end{abstract}

\maketitle

\section*{$S^1$-equivariant Floer theory}

We recall here our old (1994) heuristic construction \cite{GiZ,GiH} which highlights the role of $D$-modules in quantum cohomology theory, and adjust the construction to the case of quantum K-theory and $D_q$-modules, following the more recent exposition \cite{GiTo}.  

Let $X$ be a compact symplectic (or K\"ahler) target space, which for 
simplicity is assumed simply-connected, 
and such that $\pi_2(X)=H_2(X)\cong \ZZ^K$. Let $d=(d_1,\dots,d_k)$ be integer coordinates on $H_2(X)$, and $\omega_1, \dots, \omega_K$ be closed 2-forms on $X$ with integer periods, representing the corresponding basis of $H^2(X,\RR)$.

On the space $L_0X$ of contractible parameterized loops $S^1\to X$, as well 
as on its universal cover $\tilde{L_0X}$, one defines closed 2-forms 
$\Omega_i$, which associates to two vector fields $\xi$ and
$\eta$ along a given loop the value
\[ \Omega_i(\xi,\eta) := \oint \omega_i(\xi(t),\eta(t))\ dt.\]  
A point $\gamma \in \tilde{L_0X}$ is a loop in $X$ together with a homotopy
type of a disk $u: D^2\to X$ attached to it. One defines the 
{\em action functionals} $H_i: \tilde{L_0X}\to \RR$ by evaluating 
the 2-forms $\omega_i$ on such disks:
\[ H_i(\gamma) := \int_{D^2} u^*\omega_i .\]
  
Consider the action of $S^1$ on $\tilde{L_0X}$, defined by the rotation 
of loops, and let $V$ denote the velocity vector field of this action.
It is well-known that $V$ is
$\Omega_i$-hamiltonian with the Hamilton function $H_i$, i.e.:
\[ \operatorname{i}_{V} \Omega_i + dH_i = 0,\ \ \ i=1,\dots, K.\]

Denote by $z$ the generator of the coefficient ring $H^*(BS^1)$ of 
$S^1$-equivariant cohomology theory. The $S^1$-equivariant de Rham complex
(of $\tilde{L_0X}$ in our case) consists of $S^1$-invariant differential 
forms with coefficients in $\RR [z]$, and is equipped with the differential
$D:= d+zi_V$. Then the degree-2 elements 
\[ p_i:=\Omega_i + z H_i, \ \ \ \ i=1,\dots, K,\]
are $S^1$-equivariantly closed: $Dp_a=0$. This is standard in the context of Duistermaat--Heckman's formula.   

Furthermore, the lattice $\pi_2(X)$ acts by deck transformations on the 
universal covering $\tilde{L_0X}\to L_0X$. Namely, an element $d\in \pi_2(X)$
acts on $\gamma \in \tilde{L_0X}$ by replacing the homotopy type $[u]$ of the
disk with $[u]+d$. We denote by $Q^d=Q_1^{d_1}\cdots Q_K^{d_K}$ the operation 
of pulling-back differential forms by this deck transformation. It is an 
observation from \cite{GiZ,GiH} that the operations $Q_i$ and the operations
of exterior multiplication by $p_i$ do not commute:
\[ p_i Q_{i'} - Q_{i'} p_i = -z Q_i \delta_{ii'}.\]
These are commutation relations between generators of the algebra of 
differential operators on the $K$-dimensional torus:
\[ \left[ -z Q_i\p_{Q_i} , Q_{i'} \right] = -z Q_i \delta_{ii'}.\] 
Likewise, if $P_i$ denotes the $S^1$-equivariant line bundle on $\tilde{L_0X}$
whose Chern character is $e^{-p_i}$, then tensoring vector bundles by $P_i$ 
and pulling back vector bundles by $Q_i$ do not commute:
\[ P_i Q_{i'} = ((q-1)\delta_{ii'}+1)Q_{i'} P_i .\]
These are commutation relations in the algebra of finite-difference
operators, generated by multiplications and translations:
\[ Q_i\mapsto Q_i\times, \ \ P_i\mapsto e^{zQ_i\p_{Q_i}}=q^{Q_i\p_{Q_i}}, 
 \ \ \text{where}\ \ q=e^z.\] 
Thinking of these operations acting on $S^1$-equivariant Floer theory of the 
loop space, one arrives at the conclusion that $S^1$-equivariant Floer 
cohomology (K-theory) should carry the structure of a module over
the algebra of differential (respectively finite-difference) operators.
We will elucidate this conclusion with toric examples after giving a convenient description of toric manifolds. 

\section*{Toric manifolds}

Fans and momentum polyhedra are two the most popular languages in algebraic and symplectic geometry of toric manifolds \cite{Au}. In symplectic {\em topology}, a third framework, where toric manifolds are treated as symplectic reductions or GIT quotients of a linear space, turns out to be more convenient \cite{GiA}.

Let $\Delta$ be the momentum polyhedron of a compact symplectic toric manifold (we remind that it lives in the dual of the Lie algebra of a compact torus, and is therefore equipped with the integer lattice), and $N$ be the number of its hyperplane faces. The corresponding $N$ supporting affine linear functions with the minimal integer slopes canonically embed $\Delta$ into the first orthant $\RR_{+}^N$ in $\RR^N$, and thereby represent the toric manifold as the symplectic quotient of $\CC^N$.

Indeed, the torus $T^N$ acts by diagonal matrices on $\CC^N$ with the momentum
map $(z_1,\dots,z_N)\mapsto (|z_1|^2, \dots, |z_N|^2)$: $\CC^N \to \RR_{+}^N\subset \Lie^*T^N$. For a subtorus $T^K\subset T^N$, the momentum map is obtained by
further projection $\m: \Lie^*T^N \to \Lie^*T^K = \RR^K$. The last equality uses a basis, $(p_1,\dots, p_K)$, which we will assume integer. In fact one only needs to look at the {\em picture} $\m (\RR_{+}^N) \subset \RR^K$ of the first orthant (see example on Figure 1), i.e. to know the images $u_1,\dots, u_N$ in $\RR^K$ of the unit coordinate vectors from $\RR^N$:
\[  u_j = p_1 m_{1j}+\dots + p_K m_{Kj}, \ \ j=1,\dots, N.\]
When $\Delta$ is the fiber $ \m^{-1}(\omega)$ in the first orthant over some regular value $\omega$, the initial toric manifold is identified with the symplectic reduction $X = \CC^N//_{\omega} T^K$. Alternatively, removing from $\CC^N$ all coordinate subspaces whose moment images do not contain $\omega$,
one identifies $X$ with the quotient $\mathring{\CC}^N/T^K_{\CC}$ of the rest by the action of the compexified torus (GIT quotient), and thereby equips $X$ with a complex structure.  

\begin{figure}[htb]
\begin{center}
\epsfig{file=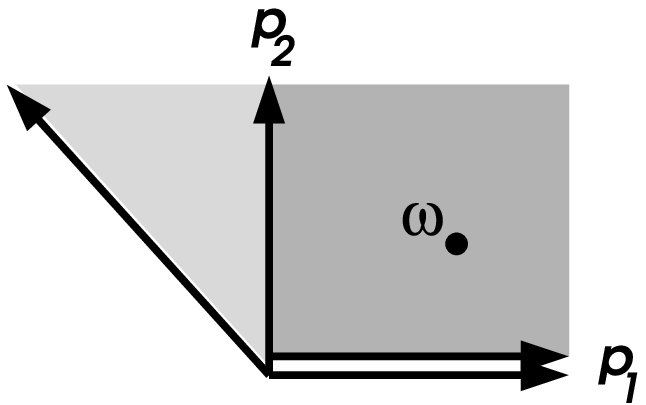} 
\end{center}
\caption{$X=\operatorname{proj}(\O_{\CC P^1}(-1) \oplus \O_{\CC P^1})$}
\end{figure}

Here is how basic topological information about $X$ can be read off the {\em picture}. The space $\RR^K$ and the lattice spanned by $p_i$ are identified with
$H^2(X, \RR)\supset H^2(X, \ZZ)$. The vectors $u_1,\dots, u_N$ represent cohomology classes of the toric divisors of complex codimension $1$ (they correspond to the hyperplane faces of the momentum polyhedron), and $c_1(T_X)$
is their sum. In the example of
Figure 1, $u_1=u_2=p_1, u_2=p_2, u_4=p_2-p_1$. The {\em chamber} (connected component) of the set of regular values of the moment map, which contains $\omega$ (it is the darkest region on Figure 1) becomes the K\"ahler cone of $X$.
It is the intersection of the images of those $K$-dimensional walls $\RR^K_{+}$
of the first orthant $\RR_{+}^N$ which contain $\omega$ in their image. In the example, there are $4$ of these: spanned by $(u_1,u_3)$, $(u_2, u_3)$, $(u_1, u_4)$,
and $(u_2, u_4)$. They are in one-to-one correspondence with the vertices of the momentum polyhedron, and hence with fixed points of $T^N$ in $X$. By the way, $X$ is non-singular if and only if the determinants of these maps $\RR^K_{+}\to \RR^K$ (i.e. appropriate $K\times K$ minors of the $K\times N$ matrix $\m$) are equal to $\pm 1$. 

The ring $H^*(X)$ is multiplicatively generated by $u_1,\dots, n_N$, which besides the $N-K$ linear relations (given by the above expressions in terms of $p_1,\dots, p_K$,) satisfy multiplicative {\em Kirwan's relations}. Namely,
$\prod_{j\in J} u_j =0$ whenever the toric divisors $u_j$ with
$j\in J\subset \{ 1,\dots, N\}$ have empty geometric intersection. In minimalist form, for each {\em maximal} subset $J\subset \{1,\dots, N\}$ such that the cone spanned by the vectors $u_j$ on the {\em picture} misses the K\"ahler cone, there is one Kirwan's relation $\prod_{j\notin J} u_j =0$. In our example, there are two Kirwan's relations: $u_1u_2=0$ and $u_3u_4=0$, i.e. the complete presentation of $H^*(X)$ is $p_1^2=p_2(p_2-p_1)=0$.

The spectrum of the algebra defined by Kirwan's relations (we call it {\em hedgehog}, after Czech, or anti-tank hedgehogs), is described geometrically as follows. For each $\alpha \in X^T$, consider the corresponding $K$-dimensional wall of $\RR_{+}^N$ whose picture contains $\omega$, and let $J(\alpha)$ denote the corresponding cardinality-$K$ subset of $\{ 1,\dots, N\}$. In the complex space with coordinates $u_1,\dots, u_N$, consider the $N-K$-dimensional coordinate subspace ({\em rail}) $\CC^{N-K}_{\alpha}$ given by the equations $u_j=0$, $j\in J(\alpha)$. The hedgehog is the union of the rails. Respectively
$H^*(X, \CC)$ is the algebra of functions on the ``thick point'', obtained by
intersecting the hedgehog with the $K$-dimensional range of the map $\CC^K \to \CC^N:  p \mapsto u = \m^t p$. In the $T^N$-equivariant version of the theory, this subspace is deformed into
\[ u_j(p)=\sum_{i=1}^K p_i m_{ij} - \la_j, \ \ j=1,\dots,N,\]
and for generic $\la$ intersects the hedgehog at isolated points corresponding
to the fixed points $\alpha X^T$.  Here $\la_j$ are the generators of the coefficient ring $H^*(BT^N)$. Finally, the operation of integration $H^*_{T^N}(X)\to H^*_{T^N}(pt)$ can be written (under some orientation convention) in the form of the residue sum over these intersection point:
\[ \int_X \phi (p.\la) = \sum_{\alpha \in X^T} \Res_{p: u(p)\in \CC^{N-K}_\alpha} \frac{\phi(p,\la)\, dp_1\wedge \cdots \wedge dp_K}{u_1(p)\dots u_N(p)}.\]
This follows from fixed point localization. 

In K-theory, let $P_i$ and $U_j$ be the ($T^N$-equivariant) line bundles whose Chern characters are $e^{-u_j}$ and $e^{-p_i}$ respectively. The ring $K^0_{T^N}(X)$
is described by Kirwan's relations
\[ \prod_{j\in J} (1-U_j)=0\ \text{whenever $\bigcap_{j\in J} u_j =\emptyset$ for corresponding toric divisors,} \]
together with the multiplicative relations
\[ U_j=\prod_{i=1}^kP_i^{m_{ij}}\gL_j^{-1}, \ \text{where $\gL_j=e^{-\la_j}$ are generators of $\operatorname{Repr} (T^N)$},\]
the coefficient ring of $T^N$-equivariant K-theory. The K-theoretic hedgehog,
defined by Kirwan's relations, lives in the complex torus $(\CC^{\times})^N$ with coordinates $U_1,\dots, U_N$, and is the union of subtori $(\CC^{\times})^{N-K}_\alpha$ given by the equations $U_j=1, j\in J(\alpha)$.
The trace operation $\tr_{T^N}: K^0_{T^N}(X)\to K^0_{T^N}(pt)$ takes on the residue form
\[ \tr_{T^N}(X;\Phi (P)) = \sum_{\alpha\in X^T} \Res_{P: U(P)\in (\CC^{\times})^{N-K}_\alpha}  \frac{\Phi (P)\, dP_1\wedge \cdots \wedge dP_K}{\prod_{j=1}^N(1-U_j(P))\, P_1\cdots P_K}.\]
This follows from Lefschetz' fixed point formula. In the example, we have:
$U_1=P_1/\gL_1$, $U_2=P_1/\gL_2$, $U_3=P_2/\gL_3$, $U_4=P_2/P_1\gL_4$. There are four intersections with the hedgehog: $U_1=U_3=1$, $U_2=U_3=1$, $U_1=U_4=1$, and $U_2=U_4=1$. The respective residues take the form:
\begin{align*} \frac{\Phi (\gL_1,\gL_3)}{(1-\gL_1/\gL_2)(1-\gL_3/\gL_1\gL_4)} +&\frac{\Phi(\gL_2,\gL_3)}{(1-\gL_2/\gL_1)(1-\gL_3/\gL_2\gL_4)}+ \\
\frac{\Phi(\gL_1,\gL_1\gL_4)}{(1-\gL_1/\gL_2)(1-\gL_1\gL_4/\gL_3)}+&
\frac{\Phi (\gL_2,\gL_2\gL4)}{(1-\gL_2/\gL_1)(1-\gL_2\gL_4/\gL_3)}.\end{align*}

\section*{Linear sigma-models}

Returning to the heuristics based on loop spaces, we replace the universal cover $\tilde{L_0X}$ of the space of contractible loops in $X=\CC^N//_{\omega}T^K$ with the infinite dimensional toric manifold $L\CC^N//_{\omega} T^K$. Note that the group $LT^K$ of {\em loops} in $T^K$ is homotopically the same as $T^K\times\pi_1(T^K)$, and that neglecting to factorize by $\pi_1(T^K)=\ZZ^K$ is equivalent to passing to the universal cover of $L_0X$. We consider the model of the loop space $L\CC^N = \CC^N [\zeta,\zeta^{-1}]$ as equivariant with respect to $T^N \times S^1$, where $T^N$ acts as before on $\CC^N$, and $S^1$ acts by rotation of the loop's parameter: $\zeta \mapsto e^{it}\zeta $. The {\em picture} in $\RR^K$ corresponding to our infinite dimensional toric manifold consists of countably many copies of each of the vectors $u_j$. The copies represent the Fourier modes of the loops, and
the equivariant classes of the corresponding toric divisors in terms of
the basis $p_1,\dots, p_K$ have the form
\[ \sum_{i=1}^K p_i m_{ij} - \la_j - rz = u_j(p)-rz, \ \ j=1,\dots,N, \ r=0,\pm 1,\pm 2,\dots\]
In K-theory of the loop space, the line corresponding line bundles are
\[ \prod_{i=1}^K P_i^{m_{ij}}\gL_j^{-1}q^{r} = U_j(P)q^{r},\ \ j=1,\dots, N,\ r=0,\pm 1,\pm 2,\dots\]
The {\em Floer fundamental cycle} $Fl_X$ in the loop space, by definition, consists of those loops which bound holomorphic disks. In our model of the loop space, $Fl_X = \CC^N[\zeta]//_{\omega}T^K$. This gives rise to the following formula for the trace over $Fl_X$:
\begin{align*} \tr_{T^N\times S^1} (Fl_X ; \Phi (P)) &= \\
\frac{1}{(2\pi i)^K}\oint \Phi(P)\ & \frac{\prod_{j=1}^N\prod_{r=1}^{\infty}(1-U_j(p)q^r)}{\prod_{j=1}^N \prod_{r=-\infty}^{\infty}(1-U_j(P)q^r)} \frac{dP_1\wedge \cdots \wedge dP_K}{P_1\cdots P_K}.\end{align*}
Thus, the structure sheaf of the semi-infinite cycle $Fl_X$ is Poincar\'e-dual to the semi-infinite product
\[ \hat{I}_X:=\prod_{j=1}^N\prod_{r=1}^{\infty}(1-U_j(P)q^r).\]
Similarly, for a toric bundle $E\to X$ or super-bundle $\Pi E$ (see Part IV),
endowed with the fiberwise scalar action of $T^1$, $\hat{I}_E$ and $\hat{I}_{\Pi E}$ are obtained from $\hat{I}_X$ by respectively division and multiplication by the K-theoretic Euler class of the obvious semi-infinite vector bundle:
\[ \hat{I}_E:=\hat{I}_X /\prod_{a=1}^L \prod_{r=0}^{\infty}(1-\la V_a(P)q^{-r}), \  \text{and} \ 
\hat{I}_{\Pi E}=\hat{I}_X\, \prod_{a=1}^L \prod_{r=0}^{\infty}(1-\la V_a(P)q^{-r}).\]
Here $\la \in T^1$, and $V_a$ are toric line bundles, 
\[ V_a(P)= P_1^{l_{1a}}\cdots P_K^{l_{Ka}}, \ \ E=\oplus_{a=1}^L V_a .\]
  
Our aim is to compute the left $D_q$-module generated by $\hat{I}_X$. The deck transformation $Q^d$ corresponding to a homology class $d\in H_2(X)$ acts in our model by $Q^d(P_i)=P_iq^{-d_i}$, where $(d_1,\dots,d_K)$ are coordinates on $H_2(X)$ in the basis dual to $(p_1,\dots, p_K)$. We find that $\hat{I}_X$ satisfies the following relations:
\[ Q^d \hat{I}_X = \prod_{j=1}^N\,  \frac{\prod_{r=-\infty}^{D_j(d)-1}(1-U_j(P)q^{-r})} {\prod_{r=-\infty}^{-1}(1-U_j(P)q^{-r})} \, \hat{I}_X^K.\]
Of course, the relations for all $Q^d$ follow from the basis relations with $Q^d=Q_i, i=1, \dots, K$. 
For instance, in our example, after some rearrangements, we obtain a system of two finite-difference equations (for $d=(1,0)$ and $(0,1)$):  
\begin{align*} (1-P_1\gL_1^{-1})(1-P_2\gL_2^{-1}) \hat{I}_X = Q_1 (1&-P_2P_1^{-1}\gL_4^{-1}) \hat{I}_X  \\
  (1-P_2\gL_2^{-1})(1-P_2P_1^{-1}\gL_4^{-1}) \hat{I}_X &= Q_2 \hat{I}_X. \end{align*} 

To save space, we refer the reader to \cite{GiH} for an explanation (though given in the cohomological context) of how to mechanically pass from this ``momentum'' representation of the Floer fundamental class (i.e. expresses as a function of $P$) to the ``coordinate'' representation in the form of
the hypergeometric $Q$-series with vector coefficients in $K^0(X)$. In that representation, $Q^d$ acts as multiplication by $Q_1^{d_1}\cdots Q_K^{d_K}$,
and $P_i$ acts as $P_iq^{Q_i\p_{Q_i}}$, i.e. as the change $Q_i\mapsto qQ_i$ 
accompanied with multiplication by $P_i$ in $K^0(X)$. With these conventions,
we have:
\[ I_X=\sum_{d\in \ZZ^K} Q^d\, \prod_{j=1}^N \frac{\prod_{r=-\infty}^0(1-U_j(P)q^r)}{\prod_{r=-\infty}^{D_j(d)}(1-U_j(P)q^r)}.\]
This is just another way to describe the same $D_q$-module, and so $I_X$
satisfies the system of finite-difference equations:
\[ \prod_{j=1}^N \, \frac{\prod_{r=-\infty}^{m_{ij}-1}(1-q^{-r}U_j(Pq^{Q\p_Q}))} {\prod_{r=-\infty}^{-1} (1-q^{-r}U_j(Pq^{Q\p_Q}))}\, I_X =  Q_i I_X, \ \ i=1,\dots, K.\]

\section*{Real life}

The toric $q$-hypergeometric function $I_X$, though comes from heuristic manipulation, has something to do with real life. 

\medskip

{\tt Theorem.} {\em The series $(1-q) I_X$ is a value of the big J-function in symmetrized $T^N$-equivariant quantum K-theory of toric manifold $X$.}

\medskip

{\tt Proof.} We follow the plan based on fixed point localization and explained in detail in Part II and Part IV in the example of complex projective spaces. 

We write $I_X=\sum_{\alpha \in X^{T^N}} I_X^{(\alpha)} \phi_{\alpha}$ is components in the basis $\{ \phi_{\alpha}\}$ of delta-functions of fixed points. Denote by $U_j(\alpha)$ the restriction of $U_j(P)$ to the fixed point $\alpha$. We have $U_j(P)=1$ for each of the $K$ values of $j\in J(\alpha)$, i.e. 
\[ P_1^{m_{1j}}\cdots P_K^{m_{Kj}}=\gL_j, \ j\in J(\alpha).\]
This determines expressions for $P_i$, and consequently for $U_j(P)$ with $j\notin J(\alpha)$ as Laurent monomials in $\gL_1,\dots, \gL_N$. We have  
\[ I^{(\alpha)}_X(q)=\sum_{d\in \ZZ_{+}^K(\alpha)} \frac{Q^d}{\prod_{j\in J(\alpha)}\prod_{r=1}^{D_j(d)}(1-q^r)}\, \prod_{j\notin J(\alpha)}\frac{\prod_{r=-\infty}^0 (1-q^rU_j(\alpha))}{\prod_{r=-\infty}^{D_j(d)}(1-q^rU_j(\alpha))}.\] 
The summation range $\ZZ_{+}^K(\alpha)$ is over $d\in \ZZ^K$ such that $D_j(d)\geq 0$ for all $j\in J(\alpha)$, because outside this range, there is a factor $(1-q^0)$ in the numerator.\footnote{It is dual in $\ZZ^K=H_2(X,\ZZ)$
to the image $\RR_{+}^K(\alpha)$ on the {\em picture} $\RR^K=H^2(X,\RR)$ of the $K$-dimensional face of $\RR_{+}^N$ corresponding to $\alpha$. Note that the intersection of all $\RR_{+}^K(\alpha)$ is the closure of the K\"ahler cone, and
respectively the convex hull of all $\ZZ_{+}^K(\alpha)$ is the {\em Mori cone}
of possible degrees of holomorphic curves.}

{\bf (i)} Temporarily encode degrees $d$ by $D_j(d), j\in J(\alpha)$, i.e. introduce Laurent monomials $Q_j(\alpha)$ in Novikov's variables $Q_1,\dots, Q_K$ such that
\[ Q_1^{d_1}\cdots Q_K^{d_K}= \prod_{j\in J(\alpha)}Q_j(\alpha)^{D_j(d)}\ \ \text{for all $d$}.\] We have:
\[ \sum_{d\in \ZZ_{+}^K(\alpha)} \frac{Q^d}{\prod_{j\in J(\alpha)}\prod_{r=1}^{D_j(d)} (1-q^r)} = e^{\textstyle \sum_{k>0} \sum_{j\in J(\alpha)}  Q_j(\alpha)^k/k(1-q^k)}.\]
According to Part I (or Part III), the right hand side is $J_{pt}(\tau)/(1-q)$
with $\tau = \sum_{j\notin J(\alpha)} Q_j(\alpha)$ in the Novikov ring considered as a $\la$-algebra with the Adams operations $\Psi^k(Q^d)=Q^{kd}$.
It follows now from results of Part IV, that $(1-q)I_X^{(\alpha)}$, expanded near the roots of unity (i.e. with the products on the right expanded as power series in $q$), represents a value of the big J-function $\J_{pt}$. Namely, recall fromPart IV that by $\Gamma$-operators we mean $q$-difference operators with symbold defined by
\[ \Gamma_{q}(x) := e^{\sum_{k>0} x^k /k(1-q^k)} \sim \prod_{r=0}^{\infty} \frac{1}{1-xq^r}\]
Expressing each $u_j$ as a linear combination $u_j=\sum_{i\in J(\a)} u_in_{ij}$, we find for multi-variable $\Gamma$-operators:
\[ \frac{\Gamma_{q^{-1}}(\la)}{\Gamma_{q^{-1}}\left(\la\, q^{\textstyle\sum_{i\in J(\a)} n_{ij}Q_i(\a)\p_{Q_i(\a)}} \right) } \, Q^d = Q^d\,
\frac{\prod_{r=-\infty}^0(1-\la q^r)}{\prod_{r=-\infty}^{D_j(d)}(1-\la q^r)}.\]
Thus, applying to $J_{pt}/(1-q)$ such operators (one for each $j=1,\dots, N$)
and setting $\la=U_j(\a)$, we obtain $I^{(\a)}_X$. Due to the invariance of the big J-function $\J_{pt}$ with respect to the $q$-difference operators (as explained in Part IV), we conclude that $(1-q) I_X^{(a)}$ is a value of $\J_{pt}$.

{\bf (ii)} All poles of $I_X^{(\alpha)}$ away from roots of unity are simple for generic values of $\gL_1,\dots, \gL_N$. We compute the residues at such poles.
The pole is specified by the choice in the denominators of one of the factors $1-q^mU_{j_0}$ with a $j_0\notin J(\alpha)$, and by the choice of one of the $m$th roots $\la^{1/m}$ of $\la:=U_{j_0}(\alpha)$. The choice of $j_0\notin J(\alpha)$ determines a 1-dimensional orbit of $T^N_{\CC}$ in $X$, connecting
the fixed point $\a$ with another fixed point, $\b$. The closure of this orbit is a holomorphic sphere $\CC P^1\subset X$, represented on the {\em picture} by
a collection of $u_j$ of cardinality $K+1$: the union 
$J(\alpha) \sqcup \{ j_0\} = J(\beta)\sqcup \{j_0'\}$, where $j_0'$ is a unique element of $J(\alpha)$ missing in $J(\b)$. The torus $T^N$ acts on the cotangent lines to this $\CC P^1$ at the fixed points by the characters $\la=U_{j_0}(\a)$ and $\la^{-1}=U_{j_0'}(\b)$ respectively (which are therefore inverse to each other). Moreover, denote by $d_{\a\b}$ the degree of this $\CC P^1$. Then $U_j(\a)/U_j(\b)=\la^{D_j(d_{\a\b})}$, and in particular $D_{j_0}(d_{\a\b})=D_{j_0'}(d_{\a\b})=1$.  Indeed, by cohomological fixed point localization on this $\CC P^1$,
\[ D_j(d_{\a\b}):=-\int_{d_{\a\b}}\ln U_j = \frac{-\ln U_j(\a)}{-\ln \la} + \frac{-\ln U_j(\b)}{\ln \la} = \frac{\ln U_j(\a)/U_j(\b)}{\ln \la}.\]
Consequently, at $q=\la^{-1/m}$ we have for all $r$ and $j$:
\[ 1-q^rU_j(\a)=1-q^{r-m D_j(d_{\a\b})}U_j(\b).\]
Under these constraints,
\[ \frac{\prod_{r=-\infty}^{mD_j(d_{\a\b})}(1-q^r U_j(\a))}{\prod_{r=-\infty}^ {D_j(d)}(1-q^rU_j(\a))} =
\frac{\prod_{r=-\infty}^0(1-q^r U_j(\b))}{\prod_{r=-\infty}^{D_j(d)-mD_j(d_{\a\b})} (1-q^rU_j(\b))} .\]
It follows that at $q=\la^{-1/m}$, 
\begin{align*} (1-q^m\la) I_X^{(\a)}(q) =(1-q^mU_{j_0}(\a))& \sum_{d\in \ZZ^K} Q^d \, \prod_{j=1}^N\frac{\prod_{r=-\infty}^{0}(1-q^rU_j(\a))}
  {\prod_{r=-\infty}^{D_j(d)}(1-q^rU_j(\a))} = \\   
 Q^{md_{\a\b}}\, (1-q^m U_{j_0}(\a)) \, &\prod_{j=1}^N \frac{\prod_{r=-\infty}^{0}(1-q^rU_j(\a))}
 {\prod_{r=-\infty}^{mD_j(d_{\a\b})} (1-q^rU_j(\a))} \times \\
 \sum_{d\in \ZZ^K}  Q^{d-md_{\a\b}} & \prod_{j=1}^N\frac{\prod_{r=-\infty}^{0}(1-q^rU_j(\b))}
     {\prod_{r=-\infty}^{D_j(d)-mD_j(d_{\a\b})}(1-q^rU_j(\b))}.     
\end{align*}
Equivalently,
\[ \Res_{q=\la^{-1/m}} I_X^{(\a)}(q)\, \frac{dq}{q}=-\frac{Q^{md_{\a\b}}}{m}\, \frac{\phi^{\a}}{C_{\a\b}(m)}\, I_X^{(\b)}(\la^{-1/m}), \]
where $\phi^{\a}=\prod_{j\notin J(\a)}(1-U_j(\a))=\Eu_{T^N}^K(T^*_{\a}X)$, and 
\[ \frac{C_{\a\b}(m)}{\phi^{\a}} = \phi^{\a} \prod_{r=1}^{m-1} (1-\la^{r/m})\, \prod_{j\neq j_0}\frac{\prod_{r=-\infty}^{mD_j(d_{\a\b})}(1-\la^{-r/m}U_j(\a))}
   {\prod_{r=-\infty}^{0} (1-\la^{-r/m}U_j(\a))} . \]
   Thus, the residues at the simple poles satisfy the recursion relations derived by fixed point localization arguments in Part II.  More precisely, it remains to check that $C_{\a\b}(m)= \Eu^K_{T^N}(T*_p\M_{0,2}(X,md_{\a,b})$, where $T*_p$ is the virtual cotangent space to the moduli space at the point $p$ represented by the $n$-multiple cover of the 1-dimensional orbit connecting fixed points $\a$ and $\b$. This verification is straightforward. In $K_{T^N}^0(X)$, we have $T^*X = U_1+\cdots+U_N-K$ (as follows from the quotient description of $X$, or by localization to $X^T$). Therefore the cotangent $T*_p$ is identified with $\oplus_j H^0(\CC P^1; U_j^m)\ominus H^1(\CC P^1;U_j^N-K-1$ (the last $-1$ stands for reparameterizations of $\CC P^1$ with $2$ marked points), which is easily described in terms of spaces of binary forms of degrees $D_j(\a\b)$. The factors in the formula for $C_{\a\b}(m)$ correspond to $T^N$-weight of the monomials in such binary forms.     

From (i) and (ii) it follows that $(1-q) I_X$ is a value of the big J-function in permutation-equivariant quantum K-theory of $X$. Since $I_X$ is defined over the $\la$-algebra $\ZZ [\gL^{\pm 1}] [[Q]]$ involving only Novikov's variables and functions on $T^N$, the value actually belongs to the {\em symmetrized} theory, i.e. carries information only about multiplicities the part of sheaf cohomology, {\em invariant} under permutations of marked points.  \qed

\pagebreak

In the case of bundle $E$ or super-bundle $\Pi E$, where $E=\oplus_{a=1}^LV_a$ is the sum of toric line bundles $V_a = \prod_i P_i^{l_{ia}}$, one similarly obtains $q$-hypergeometric series
\begin{align*} I_E&=\sum_{d\in \ZZ^K} Q^d \, \prod_{j=1}^N\frac{\prod_{r=-\infty}^0(1-q^rU_j(P))} {\prod_{r=-\infty}^{D_j(d)}(1-q^rU_j(P))}\ \prod_{a=1}^L\frac{\prod_{r=-\infty}^0(1-\la q^r V_a(P))}{\prod_{r=-\infty}^{\Delta_a(d)}(1-\la q^r V_a(P))},\\
  I_{\Pi E}&=\sum_{d\in \ZZ^K} Q^d \, \prod_{j=1}^N\frac{\prod_{r=-\infty}^0(1-q^rU_j(P))} {\prod_{r=-\infty}^{D_j(d)}(1-q^rU_j(P))}\  \prod_{a=1}^L\frac{\prod_{r=-\infty}^{\Delta_a(d)}(1-\la q^r V_a(P))}{\prod_{r=-\infty}^0(1-\la q^r V_a(P))},\end{align*}
where $\Delta_a(d)=\sum_i d_i l_{ia}$. 

\medskip

{\tt Theorem.} {\em Functions $(1-q) I_E$ and $(1-q) I_{\Pi E}$ represent some values of the big J-functions in symmetrized quantum K-theories of toric bundle space $E$ and super-bundle $\Pi E$ respectively.}

\medskip

This theorem is proved the same way as the previous one.

The above results are K-theoretic analogues of cohomological ``mirror formulas'' \cite{GiT, CGi, Ir}. The strategy we followed is due to J. Brown \cite{Br}. Some special cases were obtained in \cite{GiTo} by a different strategy. Some further results and applications can be found in the recent peprint \cite{To}.

\enddocument